\documentclass[12pt,twoside]{article}
\usepackage[latin1]{inputenc}
\usepackage{times}
\usepackage{epsfig}
\usepackage{amsfonts}
\usepackage{float}
\usepackage{amsmath}
\usepackage{amssymb}
\usepackage{latexsym}
\usepackage{graphicx, color}

\newcommand{\rnc}{\renewcommand}
\newcommand{\nc}{\newcommand}
\newcommand{\mrm}{\mathrm}
\renewcommand{\b}{\textbf}
\nc{\mb}{\mathbb}
\nc{\mac}{\mathcal}
\nc{\E}{\mb{E}}
\nc{\N}{\mb{N}}
\nc{\R}{\mb{R}}
\nc{\Q}{\mb{Q}}
\rnc{\P}{\mrm P}
\rnc{\d}{\mrm d}
\nc{\C}{\mac{C}}
\nc{\D}{\mac{D}}
\nc{\B}{\mac{B}}
\nc{\oPo}{\stackrel{\mrm P}{\rightarrow}}
\nc{\oWo}{\stackrel{w}{\rightarrow}}
\nc{\oDo}{\stackrel{\mac{D}}{\rightarrow}}
\newcommand{\leb}{\lambda \! \!  \lambda}

\relax 
\allowdisplaybreaks[3]

\relax

\newcommand{\dist}{\mbox{$\, \stackrel{d}{\longrightarrow} \,$}}
\newcommand{\prob}{\mbox{$\, \stackrel{p}{\longrightarrow} \,$}}

\newcommand{\bay}{\begin{array}}
\newcommand{\eay}{\end{array}}

\newcommand{\bqa}{\begin{eqnarray*}}
\newcommand{\eqa}{\end{eqnarray*}}

\newcommand{\bee}{\begin{eqnarray*}}
\newcommand{\eee}{\end{eqnarray*}}

\newcommand{\bea}{\begin{eqnarray*}}
\newcommand{\eea}{\end{eqnarray*}}

\newcommand{\bqan}{\begin{eqnarray}}
\newcommand{\eqan}{\end{eqnarray}}

\newcommand{\be}{\begin{eqnarray}}
\newcommand{\ee}{\end{eqnarray}}

\newcommand{\bit}{\begin{itemize}}
\newcommand{\eit}{\end{itemize}}

\newcommand{\ben}{\begin{enumerate}}
\newcommand{\een}{\end{enumerate}}

\newcommand{\beq}{\begin{equation}}
\newcommand{\eeq}{\end{equation}}

\newcommand{\bdes}{\begin{description}}
\newcommand{\edes}{\end{description}}

\newcommand{\btb}{\begin{tabular}}
\newcommand{\etb}{\end{tabular}}

\newcommand{\bcen}{\begin{center}}
\newcommand{\ecen}{\end{center}}

\newcommand{\bmp}{\begin{minipage}}
\newcommand{\emp}{\end{minipage}}

\newcommand{\vlambda}{\mbox{\boldmath $\lambda$}}

\newcommand{\vw}{\boldsymbol{w}}

\newcommand{\vN}{\boldsymbol{N}}

\newcommand{\vZ}{\boldsymbol{Z}}

\newtheorem{definition}{{\sc Definition}\sc}[section]
\newcommand{\bdefi}{\begin{definition}}
\newcommand{\edefi}{\end{definition}}

\newtheorem{appropr}[definition]{{\sc Approximation Procedure}\sc}
\newcommand{\bappr}{\begin{appropr}}
\newcommand{\eappr}{\end{appropr}}

\newtheorem{bedi}[definition]{{\sc Condition}\sc}
\newcommand{\bbd}{\begin{bedi}}
\newcommand{\ebd}{\end{bedi}}

\newtheorem{bedin}[definition]{{\sc Conditions}\sc}
\newcommand{\bbdn}{\begin{bedin}}
\newcommand{\ebdn}{\end{bedin}}

\newtheorem{corollary}[definition]{{\sc Corollary}\sc}
\newcommand{\bco}{\begin{corollary}}
\newcommand{\eco}{\end{corollary}}

\newtheorem{lemma}[definition]{{\sc Lemma}\sc}
\newcommand{\blem}{\begin{lemma}}
\newcommand{\elem}{\end{lemma}}

\newtheorem{proposition}[definition]{{\sc Proposition}\sc}
\newcommand{\bpro}{\begin{proposition}}
\newcommand{\epro}{\end{proposition}}

\newtheorem{satz}[definition]{{\sc Theorem}\sc}
\newcommand{\bsa}{\begin{satz}}
\newcommand{\esa}{\end{satz}}

\newtheorem{assumption}[definition]{{\sc Assumption}\sc}
\newcommand{\bas}{\begin{assumption}}
\newcommand{\eas}{\end{assumption}}

\newtheorem{assumptions}[definition]{{\sc Assumptions}\sc}
\newcommand{\bass}{\begin{assumptions}}
\newcommand{\eass}{\end{assumptions}}

\newtheorem{abb}{{\sc Figure}\sc}
\newcommand{\babb}{\begin{abb}}
\newcommand{\eabb}{\end{abb}}

\newenvironment{remark}{\begin{rmk}\sl}{\end{rmk}}
\newtheorem{rmk}{{\sc Remark}\sc}[section]
\newcommand{\brem}{\begin{remark}}
\newcommand{\erem}{\end{remark}}

\newenvironment{remarks}{\begin{rmks}\sl}{\end{rmks}}
\newtheorem{rmks}{{\sc Remarks}\sc}[section]
\newcommand{\brems}{\begin{remarks}}
\newcommand{\erems}{\end{remarks}}

\newenvironment{example}{\begin{exmp}\rm}{\end{exmp}}
\newtheorem{exmp}{{\sc Example}\sc}[section]
\newcommand{\bbsp}{\begin{example}}
\newcommand{\ebsp}{\end{example}}
\newcommand{\bexa}{\begin{example}}
\newcommand{\eexa}{\end{example}}

\newtheorem{model}{{\sc Model}\sc}[section]
\newcommand{\bmdl}{\begin{model}}
\newcommand{\emdl}{\end{model}}

\newtheorem{scheme}{{\sc Scheme}\sc}[section]
\newcommand{\bscm}{\begin{scheme}}
\newcommand{\escm}{\end{scheme}}

\newenvironment{tabelle}{\begin{tabl}\sl}{\end{tabl}}
\newtheorem{tabl}{{\sc Table}\sc}
\newcommand{\btab}{\begin{tabelle}}
\newcommand{\etab}{\end{tabelle}}

\newenvironment{exercise}{\begin{exc}\sl}{\end{exc}}
\newtheorem{exc}{Exercise}[section]
\newcommand{\bexe}{\begin{exercise}}
\newcommand{\eexe}{\end{exercise}}

\numberwithin{equation}{section}

\begin{document}

\title{\Large \bf \mbox{How to Bootstrap Aalen-Johansen Processes for Competing Risks?}\\ Handicaps, Solutions and Limitations.}
\author{Dennis Dobler$^*$ and Markus Pauly$^*$ 
}
\maketitle
\vfill
\vfill

\noindent${}^*$ {Heinrich-Heine University of Duesseldorf, Mathematical Institute, Germany\\

\newpage

\begin{abstract}

Statistical inference in competing risks models is often based on the famous Aalen-Johansen estimator. 
Since the corresponding limit process lacks independent increments, it is typically applied together with Lin's (1997) resampling technique involving standard normal multipliers. 
Recently, it has been seen that this approach can be interpreted as a wild bootstrap technique and that other multipliers, as e.g. centered Poissons, may lead to better finite sample performances, see Beyersmann et al. (2013). 
Since the latter is closely related to Efron's classical bootstrap, the question arises whether this or more general weighted bootstrap versions of Aalen-Johansen processes lead to valid results.
Here we analyze their asymptotic behaviour and it turns out that such weighted bootstrap versions in general possess the wrong covariance structure in the limit. 
However, we explain that the weighted bootstrap can nevertheless be applied for specific null hypotheses of interest and also discuss its limitations for statistical inference.
To this end, we introduce different consistent weighted bootstrap tests for the null hypothesis of stochastically ordered cumulative incidence functions and compare their finite sample performance in a simulation study.
\end{abstract}

\noindent{\bf Keywords:} Aalen-Johansen Estimator; Bootstrap; Competing risk; Counting processes; Cumulative incidence function; Left-truncation; Right-censoring; Weighted Bootstrap.

 \newpage

%
\section{Introduction}
%

In the widely used competing risks set-up, survival data is modeled via quite simple time continuous Markov chains, which may be described by an 
initial state (e.g. ``non-failure``) and a final state (e.g. ''failure''). Here the
latter is categorized into different absorbing states which are exclusive and may be
interpreted as the ``competing`` failure causes. 
In this context the so called cumulative incidence function (CIF), also called sub-distribution function, is of particular interest. 
For each absorbing state, i.e. failure cause, it is separately defined as the probability of occurrence for this particular failure type until a given time.
Time-simultaneous inference for the CIF is often based on its canonical Aalen-Johansen estimator, see Aalen and Johansen (1978). \nocite{AalenJoh78}
However, because of the complicated covariance structure of its standardized limit process, depending on the statistical question of interest, 
often other tools are needed to create valid statistical procedures. In this context a worthwhile and very 
promising possibility to attack this problem is the use of adequate resampling procedures like Lin's multiplier technique, see Lin (1993, 1997) \nocite{lin93}\nocite{lin97} or Martinussen and Scheike (2006) \nocite{martinussen06} for special examples with medical background.
His resampling idea is as follows: For fixed data, standard normal multipliers are introduced 
into a proper (resampling) statistic which theoretically possesses the same Gaussian limit
distribution as the corresponding normalized Aalen-Johansen process of the CIF. Then the
unknown distribution of the Aalen-Johansen process is approximated by repeatedly
generating a large number of realizations of the resampling statistic. This approach leads to
the construction of valid confidence bands, see Lin (1997). 

In the context of hypothesis testing, Bajorunaite and Klein (2007, 2008)\nocite{bajorunaite07}\nocite{bajorunaite08} as well as Sankaran et al. (2010)\nocite{sankaran10} have also studied Lin's resampling scheme to test for
equality of different CIFs in extensive simulation studies. Spitoni et al. (2012) \nocite{Spitoni_2012} investigated Lin's resampling method for estimating transition probabilities in semi-Markovian models with applications to survival analysis.

As mentioned by Cai et al (2010),\nocite{cai10} Lin's (1997) multiplier method is a special version of the general wild bootstrap approach, originally introduced by Wu (1986) \nocite{wu86} for inference in regression models. 
Recently Beyersman et al. (2013) \nocite{beyersmann12b} have provided a rigorous study of the theoretical properties of the wild bootstrap for the Aalen-Johansen estimator in competing risks allowing for independent left-truncation and right-censoring. 
There it is discussed that other multipliers such as standardized Poisson variates may help to construct more accurate confidence bands for the CIF in the competing risk set-up. 
As explained in that paper the latter is quite close in spirit to Efron's (1979) \nocite{efron79} classical bootstrap, in which the resampling scheme is generated by drawing with replacement from the sample (or an adequately transformed sample).
This motivates the question whether the classical bootstrap or other related resampling techniques may also be applied for statistical inference in one- and two-sample competing risks design. 
In particular, the current paper studies 
\begin{itemize}
 \item[(1)] the theoretical properties of a general exchangeably weighted bootstrap version of the Aalen-Johansen estimator in this context, covering amongst others the above mentioned wild bootstrap as well as Efron's original bootstrap, and
 \item[(2)] statistical applications and limitations of this general resampling approach for testing different null hypotheses of interest for the CIF. 
\end{itemize}
The weighted bootstrap approach was first introduced for i.i.d. samples by 
Mason and Newton (1992), see also \nocite{mason92} Pr{\ae}stgaard and Wellner (1993), \nocite{praestgaard93} Putter and van Zwet (1996) \nocite{Putter_1996}
as well as van der Vaart and Wellner (1996). \nocite{vaart96}
It has then been further developed and generalized to more general schemes, allowing for different dependency structures, by 
Janssen and Pauls (2003), \nocite{janssen03}
Janssen (2005), \nocite{janssen05}
del Barrio et al. (2009), \nocite{delbarrio09}
Pauly (2011) \nocite{pauly11}.
Here we focus on the technique derived in Janssen (2005) and Pauly (2011). 

Inference problems of interest in competing risk designs are given by one-, two- and $k$-sample tests for the null hypotheses of equality 
(which may correspond to the construction of time-simultaneous confidence bands) as well as of ordering of the CIF(s).
Here we focus on two-sample problems. It will turn out that for the first problem (i.e. testing equality of the CIFs of two independent groups) 
the wild bootstrap is exceptionally suited, whereas for the second problem general resampling versions of studentized Pepe (1993) tests lead to consistent inference procedures. 
The theoretical results are motivated from competing risks designs with independent left-truncation and right-censoring but will also hold 
for more general counting processes satisfying the multiplicative intensity model, see the monograph of Andersen et al. (1993) \nocite{abgk93} for more details. 

The paper is organized as follows. In Section~\ref{sec: model} we introduce the competing risks model, the CIF and its canonical Aalen-Johansen estimator. 
After recapitulating the wild bootstrap technique for these estimators, we introduce their general weighted bootstrap version in Section~\ref{sec: ResAJ} and analyze their weak convergence. 
Statistical applications for testing the null hypothesis of ordered CIFs in the two-sample case are given in Section~\ref{sec: Test} and their finite sample properties are investigated in simulations in Section~\ref{sec: sim}. 
Finally our results are discussed in Section~\ref{sec: dis} and all proofs are given in the Appendix.

%
\section{Notation, Model and Estimators}\label{sec: model}
%

To be as general as possible in the competing risks set-up we consider a
non-homogeneous Markov chain $(X_t)_{t\geq 0}$ in continuous time with finite state space~$\{0,1,\dots,k\}, k\in\mathbb{N}$. Here state~$0$ is initial with ${P}(X_0=0)=1$, and all other 
states~$1,\dots,k$, representing the competing risks, are assumed to be absorbing. For ease of convenience we restrict ourselves to the case of $k=2$ with two absorbing states. 
The corresponding transition intensities (or cause-specific hazard functions) of $(X_t)_{t\geq 0}$ from state~$0$ into state $j=1,2$ will be denoted by $\alpha_{j}$ and are assumed to exist. 
Moreover, the event time is given by ${T}=\inf\{t>0 \mid\ X_t\neq 0\}$ and allows for the following relation to the cause-specific hazards
$$
  \alpha_{j}(t) = \lim_{\Delta t \searrow 0}\frac{{P}({T}\in [t, t + \Delta t), X_{{T}}=j\,\mid\, {T} \geq t)}{\Delta t},\ \, j=1,2,
$$
with an accessible practical interpretation. 
Below we are interested in the risk development of this Markov process in time on a given interval $[0,t]$ with $t<\tau$. Here $\tau$ is a given terminal time such that 
${P}({T}> \cdot)>0$ on $[0,\tau)$ and $\tau \le \sup\{u: \int_0^u(\alpha_{1}(s)+\alpha_{2}(s)) ds < \infty \}$.
Note that the competing risk $X_{{T}}$ belongs to the set $\{1,2\}$. For exemplary practical analyses of such competing risks designs we refer the reader to Allignol et al. (2010) \nocite{Allignol_10} as well as Beyersmann et al. (2012). \nocite{beyersmann12a}

For $n$ independent replicates of this Markov chain, corresponding to the observation in time from $n$ individuals, we consider the associated bivariate counting process $\vN=(N_1,N_2)$. 
Here $N_j= \sum_{i=1}^n N_{j;i}, j=1,2,$ with
\begin{equation}
\label{eq:n}
  N_{j;i}(t)=\mathbf{1}\left(\mbox{ The $i-$th Markov chain has {observed} ''$0\mapsto j$'' transition in $[0,t]$}\right),
\end{equation}
counts the number of observed transitions into state~$j$, where $\mathbf{1}(\cdot)$ denotes the indicator function. 
It is worth to note, that, under the given assumptions, the processes $N_1$ and $N_2$ are c\`{a}dl\`{a}g and do not jump simultaneously. 
Moreover, we assume that $\vN$ fulfills the multiplicative intensity model given in Andersen et al. (1993), i.e. its 
intensity process $\vlambda=(\lambda_1,\lambda_2)$ is given by
$$
  \lambda_j =  Y(t) \alpha_{j}(t),\quad j=1,2,
$$
where $Y= \sum_{i=1}^n Y_i$ with
\begin{equation}
 \label{eq:y}
Y_i(t)= \mathbf{1}\left(\mbox{ The $i-$th Markov chain did not jump in $[0,t)$}\right)
\end{equation}
denotes the number of Markov chains without a jump shortly before time~$t$, i.e. the number at risk at $t-$. 
The assumption of a multiplicative intensity model hold, e.g., in the context of independent right-censoring, left-truncation or even filtering, see Chapter~III in Andersen et al. (1993). 
For example left-truncation means that patient~$i$ is only ``under study`` if $T_i>L_i$, i.e. its event time $T_i$ is greater than its truncation time $L_i$. 
We refer to Andersen et al. (1993) for the explicit modelling of these incomplete observations in different settings.

We are now interested in deriving statistical inference procedures for the cumulative incidence functions (CIFs), or sub-distribution functions,
$$
  F_j(t) = {P}({T}\leq t, X_{ T}=j)=\int_{0}^{t}{P}({T}>
  u-)\alpha_{j}(u)du
$$
for $j=1,\ 2.$ The corresponding sub-survival function will be denoted by $S_j(t)=1-F_j(t)$, $j=1,2$. 
Consistent estimators for the CIFs are given by the famous Aalen-Johansen estimators which are defined as
\begin{equation}
\label{eq:AJ}
\hat{F}_j(t) = \int_{0}^{t}\frac{\widehat{{P}}(T>u-)}{Y(u)} J(u) \ dN_{j}(u),
\end{equation}
for $j=1,2$. Here $J(u)=\mathbf{1}\{Y(u)>0\}$ and $\widehat{{P}}(T>u)$ denotes the Kaplan-Meier estimator. In addition, we denote the related estimator of the sub-survival function by $\hat{S}_j(t) =  1 - \hat{F}_j(t)$.
Construction of simultaneous confidence bands for a CIF, say $F_1$, are in general based on the corresponding process 
$$
  W_n(\cdot) = n^{1/2}\{\hat{F}_1(\cdot)-{F}_1(\cdot)\}
$$
which, under certain regularity assumptions, converges to a zero-mean Gaussian process. 
For example, a sufficient condition, which we will assume throughout, is the following: 
For $t < \tau$ there exists a deterministic function $y$ with $\inf_{u \in (0,t]} y(u) > 0$ such that 
\bqan
\label{eq:mainass}
 \sup_{u \in [0,t]} \left| \frac{Y(u)}{n} - y(u) \right| \prob 0.
\eqan
Here and throughout the paper, ''$\prob$'' denotes convergence in probability whereas ``$\dist$`` stands for convergence in distribution as $n \rightarrow \infty$. 
In particular, under Assumption~\eqref{eq:mainass}, the process $W_n$ inherits the following representation in terms of different local martingales
\begin{eqnarray}\label{eq:martingale rep}
  W_n(s) =  \sqrt{n} \sum_{i=1}^{n}\Big( \int_{0}^{s}\frac{S_2(u) - F_1(s)}{{Y}(u)}d{M}_{1; i}(u)
  {+ \int_{0}^{s}\frac{F_1(u) - F_1(s)}{{Y}(u)}}d{M}_{2;i}(u)\Big)+o_{P}(1),
\end{eqnarray}
where for $1\leq i\leq n, j=1,2,$
\begin{equation}
  \label{eq:martingales}
  M_{j;i}(s) = N_{j;i}(s) - \int_0^s Y_i(u) \alpha_{j}(u)\, du,
\end{equation}
are local square integrable martingales. Note, that we have suppressed the dependency on the sample size $n$ as well as the appearance of the indicator $J(u)$ in both integrals in \eqref{eq:martingale rep} for ease of convenience. 
From classical results on (local) martingales it follows from \eqref{eq:mainass} and the representation \eqref{eq:martingale rep}, see, e.g., 
Theorems~$IV.1.2$ and $IV.4.2$ in Andersen et al. (1993), that for each fixed $t<\tau$ the process $W_n$ converges in distribution on {the Skohorod space} $D[0,{t}]$ 
\begin{equation}\label{eq:weakconv}
 W_n \dist U \quad \text{on}\quad D[0,{t}]
\end{equation}
to a zero-mean Gaussian process $U$ with covariance function given by
\begin{eqnarray}\label{eq:zeta}
 \zeta({s_1},{s_2})&=&\int_0^{s_1} \frac{\{S_2(u)-F_1({s_2})\}\{S_2(u)-F_1(s_1)\} \alpha_{1}(u)}{y(u)} du\nonumber\\
   &+& \int_0^{s_1} \frac{\{F_1(u)-F_1({s_2})\}\{F_1(u)-F_1({s_1})\} \alpha_{2}(u)}{y(u)} du\label{eq:cov}
\end{eqnarray}
for $s_1\leq s_2$. 
Since the covariance function $\zeta$ is unknown and the process $U$ lacks independent increments, resampling techniques are helpful tools for developing inference procedures. 
Here Lin's resampling technique, as well as the more general wild bootstrap approach (see Beyersmann et al., 2013), attack the problem by using an adequate resampling process that in some sense 
reflects the representation \eqref{eq:martingale rep} and reproduces its distribution in the limit. This will be the starting point of the following section.
%
\section{Weighted Resampling of the Aalen-Johansen Estimator}\label{sec: ResAJ}
%
The above mentioned wild bootstrap resampling procedure aims at 
approximating the limit distribution of $W_n$ by introducing i.i.d. zero-mean random variables $G_{j;i},\,1\leq i\leq n,\,1\leq j\leq2,$ with variance $1$ and finite fourth moment into the representation \eqref{eq:martingale rep}. 
Replacing ${M}_{j; i}$ with $G_{j;i}{N}_{j; i}$ and all unknown quantities with their estimators leads to the following general wild bootstrap version of $W_n$ as introduced in Lin (1997), see also Beyersmann et al. (2013),
$$
\widehat{W}_n(s)  =  \sqrt{n} \sum_{i=1}^n \big( \int_0^s \frac{G_{1;i} (\hat{S}_2(u-) - \hat{F}_1(s))}{Y(u)} \ dN_{1;i}(u)
		     + \int_0^s \frac{ G_{2;i}(\hat{F}_1(u-) - \hat{F}_1(s))}{Y(u)} \ dN_{2;i}(u)\big),
$$
where $\hat{F}_j$ and $\hat{S}_j, j=1,2, $ are the Aalen-Johansen estimators of $F_j$ and $S_j$, respectively, see Equation \eqref{eq:AJ}. Note, that we again have suppressed the 
appearance of the indicator $J(u)$ in both integrals. 
In Beyersmann et al. (2013) it was shown that the conditional distribution of $\widehat{W}_n$ weakly converges on $D[0,{t}]$ to the same Gaussian limit process $U$
\begin{equation}\label{eq:weakwildconv}
   \widehat{W}_n \dist U \quad \text{on}\quad D[0,{t}]
\end{equation}
in probability. In practice, this result is implemented by simulating, for fixed data, a large number of independent copies of the multipliers~$G_{j;i}$, to approximate the conditional distribution of~$\widehat{W}_n$. 
Here Lin's (1997) resampling scheme is obtained for standard normal multipliers. 

To obtain a better connection with Efron's classical bootstrap we rewrite (after multiplying with $\sqrt{2}$) the above wild bootstrap statistic $\sqrt{2}\widehat{W}_n$ 
as
\begin{equation}\label{eq:wildbstat}
  \sqrt{2}\widehat{W}_n(s)  =  \sqrt{2n} \sum_{i=1}^{n} \Big(G_{1;i}X_{n;i}(s) + G_{2;i}Y_{n;i}(s)\Big) = \sqrt{2n} \sum_{i=1}^{2n} G_{i} Z_{2n;i}(s),
\end{equation}
where for $0 \leq s \leq t$ and $i = 1, \dots, n$
\bqa
  && X_{n;i}(s) = \int_0^s \frac{\hat S_2(u-)}{Y(u)} J(u) \ d N_{1;i}(u) 
    - \hat F_1(s) \int_0^s \frac{1}{Y(u)} J(u) \ d N_{1;i}(u) , \\
   && Y_{n;i}(s) = \int_0^s \frac{\hat F_1(u-)}{Y(u)} J(u) \ d N_{2;i}(u) 
    - \hat F_1(s) \int_0^s \frac{1}{Y(u)} J(u) \ d N_{2;i}(u), 
\eqa
$G_i = G_{1;i}\mathbf 1(i \leq n) + G_{2;i-n} \mathbf 1(i > n)$ and $Z_{2n;i} := X_{n;i} \mathbf 1(i \leq n) + Y_{n;i-n} \mathbf 1(i > n)$. 
Now, for fixed $s$, the representation in \eqref{eq:wildbstat} may be interpreted as a wild bootstrap version of the linear statistic $\sqrt{2n} \sum_{i=1}^{2n} Z_{2n;i}(s)$ in the array of real valued random variables 
$\vZ_{2n}(s)=(Z_{2n;i}(s))_{i\leq 2n}$. Now recall from Mammen (1992) that for linear statistics in independent observations, the consistency of the wild bootstrap and Efron's bootstrap go hand in hand. 
Translating the above representation to the classical bootstrap, where given the observations a random sample 
$Z_{2n;1}^*(s),\dots Z_{2n;2n}^*(s)$ is drawn with replacement from $\vZ_{2n}(s)$, the statistic $\widehat{W}_n^E(s) = \sqrt{n} \sum_{i=1}^{2n} (Z_{2n;i}^*(s)-\overline{Z}_{2n}(s))$ 
can be interpreted as a bootstrap version of $W_n$. Here $\overline{Z}_{2n}$ denotes the mean of $\vZ_{2n}$.
Following Mason and Newton (1992) this statistic $\widehat{W}_n^E$ can be rewritten distributionally equivalently
as 
$$
  \widehat{W}_n^E=\sqrt{2n} \sum_{i=1}^{2n} (Z_{2n;i}^*-\overline{Z}_{2n})=\sqrt{2n}\sum_{i=1}^{2n} m_{2n;i} (Z_{2n;i}-\overline{Z}_{2n}) = \sqrt{2n}\sum_{i=1}^{2n} (m_{2n;i}-1) (Z_{2n;i}-\overline{Z}_{2n}),
$$
where $(m_{2n;1}, \ldots , m_{2n;2n})$ is a multinomial-$Mult(2n, 1/2n)$-distributed random vector. 
This now motivates to study a general weighted bootstrap version of $\sqrt{2}\ W_n$, namely
\begin{equation}\label{eq: exchBSAJ}
 \widehat{W}_n^* = \widehat{W}_n^*((w_{2n;i})_i,(Z_{2n;i})_i) = \sqrt{2n} \sum_{i=1}^{2n} w_{2n;i} (Z_{2n;i} - \overline{Z}_{2n}),
\end{equation}
where $\vw_{2n} = (w_{2n;1},\dots,w_{2n;2n})$ is an exchangeable vector of random variables that is independent of $\vZ_{2n}$. For example, the choice of Efron's bootstrap weights $w_{2n;i}= m_{2n;i}-1$ delivers $\widehat{W}_n^*=\widehat{W}_n^E$. 
Following Janssen (2005) and Pauly (2011) we impose the following regularity conditions on the weights for gaining convergence of all finite dimensional distributions of the process $\widehat{W}_n^*(\cdot)$ as $n\to \infty$:
\begin{equation}
\label{G2}
n^{-1/2}\max_{1\leq i \leq 2n} |w_{2n;i} - \overline{w}_{2n}| \prob 0,
\end{equation}
\begin{equation}
\label{G3}
\frac{1}{2n}\sum_{i=1}^{2n} (w_{2n;i} - \overline{w}_{2n})^{2} \prob 1,
\end{equation}
\begin{equation}
\label{G4}
w_{2n;1} - \overline{w}_{2n} \dist  Z,
\end{equation}
where $Z$ is a random variable with $\E(Z)=0$ and $Var(Z)=1$. Moreover, it turns out that sufficient conditions for the tightness of $\widehat{W}_n^*(\cdot)$ are given by
\begin{equation}
\label{G5}
\limsup_{n \rightarrow \infty} E[(w_{2n;1} - \overline{w}_{2n})^4] < \infty,
\end{equation}
\begin{equation}
\label{G6}
\limsup_{n \rightarrow \infty} n E[(w_{2n;1} - \overline{w}_{2n})^2 (w_{2n;2} - \overline{w}_{2n}) (w_{2n;3} - \overline{w}_{2n})] < \infty,
\end{equation}
\begin{equation}
\label{G7}
\limsup_{n \rightarrow \infty} n^2 E[(w_{2n;1} - \overline{w}_{2n}) (w_{2n;2} - \overline{w}_{2n}) (w_{2n;3} - \overline{w}_{2n}) (w_{2n;4} - \overline{w}_{2n})] < \infty.
\end{equation}
Heuristically, the additional Assumptions \eqref{G5}--\eqref{G7} ensure that the correlation between multiple factors of centered weights decreases quickly enough for large $n$ and a high number of different leading terms.
Under these assumptions we can prove the following weak convergence result for the exchangeably weighted bootstrap version \eqref{eq: exchBSAJ} of the Aalen-Johansen estimator.

\begin{satz}\label{theo:exchCLT}
Suppose that \eqref{eq:mainass} holds and that the weights fulfill the Conditions \eqref{G2}--\eqref{G7}. Then, for every $t<\tau$, 
the conditional distribution of $\widehat{W}_n^*(\cdot)$ given the data weakly converges on $D[0,{t}]$ to a zero-mean Gaussian process $U^*$
\begin{equation}\label{eq:weakexchnBSconv}
   \widehat{W}_n^* \dist U^* \quad \text{on}\quad D[0,{t}]
\end{equation}
in probability, where the covariance function of $U^*$ is given by $(r,s) \mapsto \zeta^*(r,s) = 2 \zeta(r,s) - \xi(r) \xi(s)$ with $\zeta$ as in \eqref{eq:zeta} and
 \begin{equation}
  \xi(s) = \int_0^s \{ S_2(u) - F_1(s) \} \alpha_{1}(u)  \ du
  + \int_0^s \{ F_1(u) - F_1(s) \} \alpha_{2}(u) \ du.
 \end{equation}
\end{satz}

\begin{remark}\label{remark:theo}\text{ }\\
\noindent(a) Note, that by means of partial integration the covariance perturbation functions $\xi$ can be rewritten as 
$\xi (s) =\int_0^s (1-A_1 -A_2) dF_1,$ where $A_j(s) = \int_0^s \alpha_j(u) \d u$ for $j=1,2$.\\[.5ex]
\noindent(b) Examples for resampling weights that fulfill the Conditions~\eqref{G2} -- \eqref{G4} can be found in Janssen (2005) as well as Pauly (2009, 2011).\nocite{pauly09} 
Weights that additionally fulfill \eqref{G5} -- \eqref{G7} can be found in Example~\ref{ex.weights} in the Appendix. 
As special examples Efron's bootstrap, the i.i.d. weighted bootstrap, as well as the Bayesian bootstrap, 
the Poisson bootstrap or even row-wise i.i.d. wild bootstrap weights $w_{2n;i}$ (with $\E(w_{2n;1})=0$, $Var(w_{2n;1})=1$ and $\limsup_{n \rightarrow  \infty} \E(w_{2n;1}^4) < \infty$) 
fulfill the limit theorem \eqref{eq:weakexchnBSconv} provided that $w_{2n;1} \dist Z$.\\[.5ex]
\noindent(c) The above theorem shows that the weighted bootstrap with exchangeable weights leads to a bootstrap version of $W_n$ 
  whose limit covariance function differs from the \emph{correct} asymptotical covariance of the Aalen-Johansen process $W_n$ by the summand $\frac{1}{2}\xi(r) \xi(s)$.\\[.5ex]
\noindent(d) In comparison, the wild bootstrap statistic $\widehat{W}_n$ from the beginning of Section~\ref{sec: ResAJ} reproduces the correct limit process. The reason for this behaviour can easily be explained at the special case of the classical bootstrap version 
(and also holds for many other related resampling versions that fall into our approach). Efron's bootstrap version of a linear statistic namely needs the involved centering of each random variable $Z_{2n;i}$ at the mean $\overline{Z}_{2n}$. 
Without this term, the bootstrap statistic $\sqrt{2n} \sum_{i=1}^{2n} m_{2n;i} Z_{2n;i}$ (with conditional expectation $(2n)^{3/2} \overline{Z}_{2n}$) would in general not follow a non-degenerated conditional limit theorem. 
However, this centering affects the (conditional) covariance structure of the bootstrap process. In particular, it can be seen in the appendix, that its asymptotic covariance function $\zeta^*(r,s)$ is given by the 
limit (in probability) of $\sum_{i=1}^{2n} 2n (Z_{2n;i}(r) - \bar Z_{2n}(r))(Z_{2n;i}(s) - \bar Z_{2n}(s))$. 
In comparison the asymptotic covariance function of the wild bootstrap version $\sqrt{2} \widehat{W}_n$ is given by 
the limit (in probability) of $\sum_{i=1}^{2n} 2n Z_{2n;i}(r)Z_{2n;i}(s)$, see the proof of Theorem~2 in Beyersmann et al. (2013).
The reason is that due to the i.i.d. structure of the zero-mean wild bootstrap weights no centering term is needed to gain a conditional central limit theorem. 
Actually, Theorem~\ref{theo:exchCLT} even shows that a resampling version of the Aalen-Johansen estimator of the form \eqref{eq: exchBSAJ} with a sequence of i.i.d. wild bootstrap weights $(w_{2n;i})_i$ 
would not possess the correct limit structure due to involved centering term $\overline{Z}_{2n}$.

\end{remark}

This result now leads to the question whether Efron's bootstrap (or other included resampling techniques that fall into our approach) is not applicable for statistical inference about CIFs in competing risks studies. 
The answer is two-fold. Since $\widehat{W}_n^*$ reproduces the wrong covariance of the Aalen-Johansen estimator it is not applicable directly. 
This means that the asymptotic limit distribution of transformed versions (as $\sup$-distances or integral statistics) of the Aalen-Johansen estimator of a CIF that serve as test statistic for a particular problem (as testing equality or ordering of a CIF) 
can in general not be reproduced by its corresponding transformed exchangeably weighted bootstrap version \eqref{eq: exchBSAJ}. 
However, for some situations it may nevertheless be applicable by including adequate studentizations to the corresponding test statistic, see e.g. Janssen (1997) \nocite{janssen97} or Pauly et al. (2012) \nocite{pauly12} for similar examples in the context of testing. 
Such a multplicative studentization works, e.g., if the statistic we are interested in becomes asymptotically pivotal after studentizing. 

To explain this statement we give a negative and a positive example. First, let us exemplify Cram\'{e}r-van Mises-type statistics for testing equality of a CIF. 
In this case the asymptotic limit is given by a squared $L_2$-norm of a Gaussian process which admits a principal components decomposition and its covariance function is a series 
depending on all eigenfunctions and eigenvalues of a corresponding integral operator, see Adler (1990) \nocite{Adler90} or Shorack and Wellner (2009) \nocite{shorack09} for details. 
In this case it seems reasonable that one studentization alone cannot transform this random variable into another principal components decomposition with predefined eigenvalues and eigenfunctions. 
Hence the result from Theorem~\ref{theo:exchCLT} is not applicable in this situation. 
However, if we consider, e.g., a one- or two-sample version of Pepe's test for the hypothesis of ordered CIFs, then it turns out that the resulting test statistic is asymptotically normal. 
In this situation a studentized version of the test statistic leads to an asymptotic standard normal distribution (in the non-degenerated case) and its finite sample distribution may be approximated by a related 
studentized bootstrap version. This will be studied in more detail in the next section for the more interesting two-sample case.

%
\section{Two-Sample Resampling Tests for Ordered CIFs}\label{sec: Test}
%
In order to demonstrate the applicability of the above theory we study a specific inference problem of interest. Suppose we are interested in the comparison of two CIFs on a subinterval $[t_1, t_2]$ of $[0,\tau)$ with $0\leq t_1< t_2<\tau$. 
Here we like to test whether the CIFs from two independent groups with the same competing risk, say $j=1$, possess a specific order. 
A practical interpretation may be given by two independent medical studies for the side effects of similar but different drugs. Another example is given in Bajorunaite and Klein (Example~$5$, 2007) where 
bone marrow transplant studies are compared. Note that similar null hypotheses (mainly the null hypothesis of equality) have already been studied in the literature, see e.g. 
Gray (1988), \nocite{Gray_1988} Aly et al. (1994), \nocite{emad_1994} Barmy et al. (2006), \nocite{barmy_2006} Bajorunaite and Klein (2007, 2008) or Sankaran et al. (2010) and the references cited therein, where some of them also apply Lin's resampling technique.

In the sequel we extend the notation from Section~\ref{sec: model} with a superscript ${}^{(k)}$ to denote the quantities of the $k$th group, $k=1,2$. This yields the CIFs 
$F_1^{(k)}$ for the competing risk $j=1$ 
as well as counting processes 
$N_{j;i}^{(k)}, Y_i^{(k)}, 1\leq i\leq n_k$, where $n_k$ is the sample size of group $k=1,2$. 
The hypotheses of interest may than be written as
\begin{eqnarray*}
 H_\leq: \{ F^{(1)}_{1} \leq F^{(2)}_{1} \text{ on } [t_1,t_2] \} 
 \quad \text{versus} \quad 
 K_\gneqq: \{ F^{(1)}_{1} \gneqq F^{(2)}_{1} \text{ on } [t_1,t_2]  \}.
\end{eqnarray*}
To this end, we suggest an integral-type test statistic, namely
\begin{equation}\label{eq: tn}
 T_n= \int_{t_1}^{t_2} \rho(u) \sqrt{ \frac{n_1 n_2}{n} }(\hat F^{(1)}_{1}(u) - \hat F^{(2)}_{1}(u)) \d u,
\end{equation}
where $n=n_1+n_2$ and $\rho: [0,\tau] \rightarrow (0,\infty)$ is a deterministic and integrable function that allows for different weighting of time intervals of interest, see e.g. Pepe (1993) for a similar choice. 
Note, that such statistics are motivated from related goodness of fit problems, see, e.g., Shorack and Wellner (2009) or van der Vaart and Wellner (1996). \nocite{shorack09}
Well known theorems from stochastic process theory then show that $T_n$ is asymptotically $N(0,\sigma_\zeta^2)$-distributed  under $\{F_1^{(1)} =F_1^{(2)}\} $ provided that $n_k/n\rightarrow p_k\in(0,1)$ for $k=1,2$.
Here the limit variance is given by 
\begin{equation}\label{eq: varT}
 \sigma_\zeta^2= \int_{t_1}^{t_2} \int_{t_1}^{t_2} \rho(s) (p_2\zeta^{(1)} + p_1\zeta^{(2)}) (s,t) \rho(t) \d s \d t,
\end{equation}
where $\zeta^{(k)}$ denotes the asymptotic covariance function of the Aalen-Johansen process $W_{n_k}^{(k)}$ of group $k=1,2$, see Equation~\eqref{eq:zeta} above. 
Note, that $\sigma_\zeta^2>0$ holds if we have $\alpha_1^{(k)}>0$ on a set with positive Lebesgue-$\leb_{|[t_1,t_2]}$ measure for at least one choice of $k=1,2,$ which we like to assume in the sequel. 
As already explained at the end of Section~\ref{sec: ResAJ} we need an asymptotically pivotal test statistic for applying our weighted bootstrap result from Theorem~\ref{theo:exchCLT}. 
This will be done by studentizing $T_n$ and will correct for the wrong bootstrap limit covariance. To this end, we construct a consistent estimate $V_{n}^2$ by replacing 
$p_2\zeta^{(1)} + p_1\zeta^{(2)}$ in \eqref{eq: varT} with $\hat \zeta_{n} := \frac{n_2}{n} \hat \zeta_{n_1}^{(1)} + \frac{n_1}{n} \hat \zeta_{n_2}^{(2)}$. 
Thereby $\hat \zeta_{n_k}^{(k)}$ is constructed by substituting the unknown CIFs $F_j^{(k)}(u)$, intensities $\alpha_j^{(k)}(u)\d u$ and the function $y^{(k)}$ in $\zeta^{(k)}$ with their canonical estimators 
$\hat F_j^{(k)}(u -)$, $\d \hat A_j^{(k)}(u)$ (the increments of the Nelson-Aalen estimator) and $Y^{(k)}/n_k$. Then, as shown in Theorem~\ref{theo: tests} below, an asymptotic level $\alpha$ test is given by 
$$
  \varphi_n=\mathbf{1}\{T_{n,stud} >u_{1-\alpha}\},
$$ 
where $u_{1-\alpha}$ denotes the $(1-\alpha)$-quantile of the standard normal distribution and $T_{n,stud}=T_n/V_{n}\mathbf{1}\{V_n>0\}$. We will now construct a weighted resampling version of $\varphi_n$. 
In view of Theorem~\ref{theo:exchCLT} and the martingale representation \eqref{eq:martingale rep} under $\{ F_1^{(1)} = F_1^{(2)}\} $ a weighted resampling version of $T_n$ may be given by
\begin{equation}\label{eq: tn*}
 T_n^* = \sqrt{\frac{n_1n_2}{n}} \sum_{k=1}^2 \sum_{i=1}^{2n_k} \int_{t_1}^{t_2} \rho(u) w_{2n;i}^{(k)} (Z_{2n;i}^{(k)}(u) - \overline{Z}_{2n}(u)) \d u,
\end{equation}
where $(w_{2n;i}^{(k)})_{i,k}$ is an array of exchangeable weights fulfilling \eqref{G2} -- \eqref{G7} and we set
$\overline{Z}_{2n} = \frac{1}{2n}  \sum_{k=1}^2 \sum_{i=1}^{2n_k}Z_{2n;i}^{(k)}$ with $Z_{2n;i}^{(k)}= (-1)^{k+1}(X_{n_k;i}^{(k)}\mathbf{1}(i\leq n_k) + Y_{n_k;i-n_k}^{(k)}\mathbf{1}(n_k< i))$. 
We like to note, that the $(-1)$ in this expression is due to the martingale representation of $T_n$. 
As shown below, an application of Theorem~\ref{theo:exchCLT} yields that the conditional distribution of $T_n^*$ is asymptotically $N(0,\sigma_{\tilde{\zeta}}^2)$-distributed in probability, where 
$\sigma_{\tilde{\zeta}}^2\neq \sigma_\zeta^2$ due to the wrong limit covariance structure of the weighted bootstrap Aalen-Johansen estimator. \\
As has already been seen in Janssen (2005) as well as Konietschke and Pauly (2013)\nocite{kopa13}, different, say classes, of weights need different studentizations. 
For ease of convenience, and to avoid distinguishing between too many cases, we therefore now focus only on two resampling procedures: Efron's bootstrap with weights $w_{2n;i}= m_{2n;i}-1$ and the wild bootstrap with 
$w_{2n;i}=G_i$. Here $(m_{2n;1},\dots,m_{2n;2n})$ is a multinomially distributed random vector with sample size $2n=\sum_{i=1}^{2n} m_{2n;i}$ and equal selection probability $1/2n$ and 
$(G_i)_i$ is a sequence of i.i.d. random variables with $\E(G_1)=0, Var(G_1)=1$ and $\E(G_1^4)<\infty$. However, other resampling tests can be obtained similarly. 
Motivated from the weighted variance estimator given in Janssen (2005, Section~3), a weighted resampling version of $V_{n}^2$, say $V_n^{*\, 2}$, is then given 
by replacing $p_2\zeta^{(1)} + p_1\zeta^{(2)}$ in \eqref{eq: varT} with $\zeta_n^* - \xi_n^*$, where
\begin{eqnarray*}
 \zeta_n^*(s,t) &=&  \frac{n_1 n_2}{n} \sum_{k=1}^2 \sum_{i=1}^{2n_k} v_{2n;i}^{(k)} Z_{2n;i}^{(k)}(s) Z_{2n;i}^{(k)}(t),\\
 \xi_n^*(s,t) &=& \frac{n_1 n_2}{2 n^2}\Big( \sum_{k=1}^2 \sum_{i=1}^{2n_k} v_{2n;i}^{(k)} Z_{2n;i}^{(k)}(s) \Big) 
    \Big( \sum_{k=1}^2 \sum_{i=1}^{2n_k} v_{2n;i}^{(k)} Z_{2n;i}^{(k)}(t) \Big).
\end{eqnarray*}
We thereby choose $v_{2n;i}=m_{2n;i}$ in case of Efron's and $v_{2n;i}=G_i^2$ in case of the wild bootstrap.
With this choice it is proven in the appendix that, under $H_= :\{F_1^{(1)}=F_2^{(2)} \text{ on } [t_1,t_2]\}$ and the conditions given in Theorem~\ref{theo: tests} below, 
the conditional distribution of $T_{n,stud}^*=T_n^*/V_n^* \mathbf{1}\{V_n^*>0\}$ given the data is asymptotically $N(0,1)$-distributed in probability. Moreover, the resulting weighted resampling tests (corresponding either to Efron's or wild bootstrap weights)
$$
  \varphi_n^*=\mathbf{1}\{T_{n,stud} > c_n^*(\alpha)\},
$$ 
are consistent and even asymptotically effective, where $c_n^*(\alpha)$ is the (data-dependent) $(1-\alpha)$-quantile of the conditional distribution of $T_{n,stud}^*$ given the data.

\begin{satz}\label{theo: tests}
Suppose that \eqref{eq:mainass} holds for both groups. Then $\varphi_n$ is a consistent and asymptotic level $\alpha$ test, i.e. $E_{H_\leq}(\varphi_n)\rightarrow \alpha\mathbf{1}\{F_1^{(1)}=F_2^{(2)}\}$ and 
$E_{K_\gneqq}(\varphi_n) \rightarrow 1$. If in addition $\sigma_{\tilde{\zeta}}^2>0$ then $\varphi_n^{*}$ is also consistent and of asymptotic level $\alpha$. 
Moreover, $\varphi_n$ and $\varphi_n^{*}$ are even asymptotically equivalent, i.e. under $H_= \,$ it holds 
$
 E_{H_=}(|\varphi_n-\varphi_n^{*}|)\rightarrow 0.
$
\end{satz}

\begin{remark}\label{Rem: Tests} (a) The asymptotic equivalence implies that both tests also possess the same power under contiguous alternatives.\\[0.5ex]
(b) In case of the wild bootstrap the results remain valid if we omit the centering term $\overline{Z}_{2n}$ in \eqref{eq: tn*} as well as the covariance correction $\xi_n^*(s,t)$. 
Below we will denote the resulting test as $\varphi_n^W$.\\[0.5ex]
(c) Note that the assumption of a deterministic weight function can be relaxed. In particular, it can be shown that the above theorem remains also valid for non-deterministic sequences of weights $\rho_n:[0,\tau]\rightarrow (0,\infty)$ such that 
$\sup_s|\rho_n(s)-\rho(s)|\oPo 0$ in probability for an integrable and deterministic function $\rho:[0,\tau]\rightarrow (0,\infty)$. 
This can be shown using straightforward stochastic process arguments similar to those applied in Brendel et al. (2013).\nocite{Brendel_2013}\\[0.5ex]
(d) Utilizing the squared weights $v_{2n;i}=G_i^2$ within the wild bootstrap variance estimator can be motivated from corresponding symmetry-type tests with weights $G_i = \frac{1}{2}(\varepsilon_1 + \varepsilon_{-1})$. 
Such tests are typically applied in the context of paired data, where the involved studentization of the test statistic is often invariant under reflections of the coordinates, see Janssen (1999) or Konietschke and Pauly (2013) for details and examples. \nocite{janssen99}
In this case, the resampling (symmetry-type) version of the studentization remains unchanged, which here corresponds to the case $G_i^2=1$ for this choice of weights. 
Hence the choice with $v_{2n;i}=G_i^2$ generalizes this to all covered wild bootstrap procedures.
\end{remark}
In the next section the finite sample properties of the asymptotic test $\varphi_n$, Efron's bootstrap test $\varphi_n^E$ ($=\varphi_n^*$ with weights $w_{2n;i}= m_{2n;i}-1$) and 
the Wild bootstrap test $\varphi_n^W$ from Remark~\ref{Rem: Tests} with normal multipliers are investigated in a small Monte Carlo study.

%
\section{Simulations}\label{sec: sim}
%

The testing procedures from the last section are all valid asymptotically, i.e. as $n\rightarrow \infty$. 
In the next step their small sample properties are investigated in a small simulation study with regard to 
(i) keeping the preassigned error level under the null hypothesis and (ii) to their power behaviour under certain alternatives. All simulations were conducted with the help of the R-computing environment, version 2.15.0 (R Development Core Team, 2010), each with
$N_{sim} = 1000$ simulation runs. Moreover, for the resampling tests we have additionally run $B= 999$ bootstrap runs in each simulation step.
Here we consider the following simulation set-up for the type-I-error:
\begin{enumerate}
\item For the event times we have modeled the cause specific intensities  of the first group as $\alpha_{1}^{(1)}(u) = \exp(-u)$ and $\alpha_{2}^{(1)}(u) = (1-\exp(-u))$ and for the second group as 
$\alpha_{1}^{(2)} \equiv c \equiv 2-\alpha_{2}^{(2)}$, where $0\leq c\leq1$ holds. Here the case $c=1$ corresponds to the situation under the null with equal CIFs of the first risk and $c<1$ implicates the alternative.
\item As sample sizes we have chosen $(n_1,n_2)= (50,50), (50,100), (100, 100)$ and let $[t_1,t_2] = [0,1.5]$ be the domain of interest.
\item Moreover, each setting was simulated both with and without right-censoring, where the censoring times were simulated as independent exponential random variables with parameter $\lambda^{(k)}$ and 
pdf $f^{(k)}(x) = \lambda^{(k)} \exp(-\lambda^{(k)} x) \mathbf{1}_{(0,\infty)}(x)$  in group $k$. 
In case of censoring we have analyzed situations with equal censoring $(\lambda^{(1)},\lambda^{(2)})= (0.5,0.5)$ (light censoring) and $(\lambda^{(1)},\lambda^{(2)})= (1,1)$ (moderate censoring) as well as unequal censoring distributions 
with $(\lambda^{(1)},\lambda^{(2)})= (0.5,1)$.
\end{enumerate}

The results for the type I errors (for $\alpha = 0.05$) of the three tests can be found in Table~\ref{table: sim}, where 
the case without censoring is denoted by $(\lambda_1,\lambda_2)=(0,0)$. For easier reading the closest result to the prescribed $5\%$ level is printed in bold type.
Note, that in this setting we have equality of the CIFs $F_1^{(k)}(t)= 0.5(1-\exp(-2t)), k=1,2,$ of the first risk $j=1$ but unequal CIFs of the second risk.
It is seen that, for most of the scenarios, the bootstrap test $\varphi_n^E$ based on Efron's multinomially distributed weights
has a simulated type I error far above the $5\%$ level (sizes in $[.049,.074]$).
Thus, $\varphi_n^E$ tends to be quite liberal.
On the contrary, the test $\varphi_n$ based on the $95\%$-quantile of the standard normal distribution, 
and the wild bootstrap test $\varphi_n^W$ based on i.i.d. standard normally distributed weights keep the $5\%$ level much better. 
In most cases, $\varphi_n^W$ (sizes in $[.041, .062]$) seems to be slightly more accurate than $\varphi_n$ (sizes in $[.041, .063]$), 
especially in settings with unbalanced sample sizes $(n_1,n_2) = (50,100)$.

\begin{table}
\begin{center}
\setlength{\tabcolsep}{5pt}
\begin{tabular}[h]{cc|c|c|c|c|c|c|c|c|c}
\hline
 & \multicolumn{1}{r|}{$(n_1,n_2)$} & \multicolumn{3}{|c|}{(50,50)} & \multicolumn{3}{|c|}{(50,100)} & \multicolumn{3}{|c}{(100,100)} \\
\multicolumn{2}{c|}{$(\lambda_1,\lambda_2)$}  & $\varphi_n$ & $\varphi_n^W$ & $\varphi_n^E$ & $\varphi_n$ & $\varphi_n^W$ & $\varphi_n^E$ & $\varphi_n$ & $\varphi_n^W$ & $\varphi_n^E$ \\\hline
\multicolumn{2}{c|}{(0,0)} & .054 &{\bf .053} & .068 & .041 & .043 & {\bf.050} & .043 & .048 & {\bf.049} \\
\multicolumn{2}{c|}{(0.5,0.5)} & .045 &{\bf.048} & .056 & {\bf.060} & {\bf.060} & .069 & {\bf.051} & .054 & .062 \\
\multicolumn{2}{c|}{(0.5,1)} & .056 &{\bf.053} & .062 & .057 & {\bf.055} & .064 & {\bf.054} & {\bf.054} & .060 \\
\multicolumn{2}{c|}{(1,0.5)} & .042 &.041 & {\bf.051} & .060 & {\bf.056} & .074 & .055 & {\bf.054} & .059\\
\multicolumn{2}{c|}{(1,1)} & {\bf.053} &.054 & .063 & .063 & {\bf.062} & .072 & {\bf.054} & .056 & .062 
\end{tabular}
\caption{Simulated size of $\varphi_n$ and the resampling tests $\varphi_n^W, \varphi_n^E$ for nominal size $\alpha=5\%$ under different sample sizes and censoring distributions}
\label{table: sim}
\end{center}
\end{table}

The results for the power of all tests are presented in Table~\ref{table: power}, 
where simulations have been performed for alternative hypotheses corresponding to $c=0.1, 0.2, \dots, 0.9$. Here the choice $c=0.9$ corresponds to a situation close to the null, whereas we move farther into the alternative with decreasing $c$. 
Apparently, $\varphi_n^E$ has the greatest power in all scenarios due to its quite liberal behaviour. 
Therefore, we turn our attention to the differences in the results for $\varphi_n$ and $\varphi_n^W$.
Apart from a few exceptions, $\varphi_n^W$ has a marginal greater power than $\varphi_n$. In particular, all of the differences in the simulated powers of these two tests amount values in the interval $[-.006,.0.012]$. 

\begin{table}
\begin{center}
\setlength{\tabcolsep}{5pt}
\begin{tabular}[h]{c|c|c|c|c|c|c|c|c|c|c|c|c}
\hline
 \multicolumn{1}{r|}{$(n_1,n_2)$} & \multicolumn{6}{|c|}{(50,50)} & \multicolumn{6}{|c}{(100,100)}  \\
 \multicolumn{1}{r|}{$(\lambda_1,\lambda_2)$}  & \multicolumn{3}{|c|}{(0,0)} & \multicolumn{3}{|c|}{(1,1)} & \multicolumn{3}{|c|}{(0,0)} & \multicolumn{3}{|c}{(1,1)}\\
c  & $\varphi_n$ & $\varphi_n^W$ & $\varphi_n^E$ & $\varphi_n$ & $\varphi_n^W$ & $\varphi_n^E$ & $\varphi_n$ & $\varphi_n^W$ & $\varphi_n^E$ & $\varphi_n$ & $\varphi_n^W$ & $\varphi_n^E$\\\hline
0.9 & .121 & .127 & .142 & .106 & .111 & .133 & .163 & .167 & .171 & .126 & .134 & .146 \\
0.8 & .244 & .245 & .280 & .206 & .210 & .241 & .345 & .349 & .373 & .302 & .300 & .330 \\
0.7 &  .404 & .409 & .448 & .341 & .335 & .385 & .595 & .596 & .613 & .518 & .530 & .561 \\
0.6 &  .588 & .595 & .625 & .511 & .510 & .557 & .828 & .832 & .851 & .744 & .742 & .768 \\
0.5 &  .774 & .775 & .814 & .662 & .667 & .711 & .962 & .963 & .968 & .893 & .892 & .911 \\
0.4 &  .920 & .921 & .932 & .817 & .817 & .844 & .992 & .991 & .993 & .978 & .977 & .983 \\
0.3 &  .982 & .982 & .985 & .931 & .932 & .948 & 1 & .999 & 1 & .995 & .996 & .998 \\
0.2 &  1 & .999 & 1 & .980 & .981 & .985 & 1 & 1 & 1 & 1 & 1 & 1 \\
0.1 &  1 & 1 & 1 & .997 & .997 & .997 & 1 & 1 & 1 & 1 & 1 & 1 \\
\end{tabular}
\caption{Simulated size of $\varphi_n$ and the resampling tests $\varphi_n^W, \varphi_n^E$ for nominal size $\alpha=5\%$ under different sample sizes and censoring distributions}
\label{table: power}
\end{center}
\end{table}

Thus, having the simulated type I error rates in mind,
there is a clear preference for $\varphi_n^W$. However, since the improvement compared to $\varphi_n$ is not very large, we plan to study the behaviour of the presented tests in a more applied paper in the future, where they will be additionally 
compared with other existing procedures. 
There, also other resampling versions that fall into our approach (such as the i.i.d. weighted bootstrap, Rubin's Bayesian bootstrap 
or simply other i.i.d. weigths with finite fourth moment, cf. Example~\ref{ex.weights}) shall be studied in extensive simulations for different settings. 
On the other hand, the simulation results for the present set-up strongly suggest not to use $\varphi_n^E$ in this context.

%
\section{Discussion and Outlook}\label{sec: dis}
%

We have considered a weighted bootstrap approach for the Aalen-Johansen estimator (AJE) of a competing risk including amongst others Efron's classical, Rubin's Bayesian as well as the wild bootstrap. 
It turned out that the asymptotic covariance structure of the AJE is not reflected correctly by the weighted bootstrap. 
This handicap is due to the utilized resampling from centered data which is a necessity for most of the presented bootstrap procedures. 
One exception is the wild bootstrap of Lin (1997) and Beyersmann et al. (2013), where this centering is not needed due to the i.i.d. structure of the weights. 
Nevertheless, we have demonstrated that the covariance problem can be solved for specific inference problems. 
Roughly speaking, the general weighted bootstrap approach can be used for test statistics (here functionals of AJEs) which are asymptotic pivots. 
This has been exemplified for the unpaired two-sample testing problem of ordered CIFs. There an integral-type statistic is made asymptotically pivotal by an adequate studentization. 
If, however, the limit distribution of the test statistic is more complicated (e.g. if a variance stabilizing transformation or studentization cannot deduce pivotality), the general weighted bootstrap is 
not applicable. In such cases as, e.g., nonparametrically testing for equality of different CIFs, the (general) wild bootstrap from uncentered observations ${\bf Z}$ seems to be the only known and reasonable choice. 
To this end, other possibilities for testing equality of different CIFs than the wild bootstrap will be studied by the authors in a forthcoming paper.\\
Finally, we like to note that in semiparametric models the above approach may be improved by modifying the presented resampling algorithms as in Lin et al. (2000) \nocite{lin00} or Scheike and Zhang (2003) \nocite{scheike03}, where 
the martingale increments $dM_{0j;i}$ in the resampling step are replaced with estimated increments $d\widehat{M}_{0j;i}$ rather than $dN_{0j;i}$.

\section*{Acknowledgements}

The authors like to thank Arthur Allignol, Jan Beyersmann and Arnold Janssen for helpful discussions and Marc Ditzhaus for computational help. 
Moreover, both authors appreciate the support received by the SFF grant F-2012/375-12.

%
\section{Appendix}
%

{\it Proof of Theorem~\ref{theo:exchCLT}.} 
In order to prove the result we have to show (conditional) weak convergence of finite dimensional distributions as well as tightness. For the first we will apply Theorem~4.1 in Pauly (2011) and for the latter we 
use a tightness criterion by Billingsley (1999). \nocite{billingsley99}
To verify the finite dimensional convergence of the process let $t_1, \dots, t_k \in [0,t]$. Then, as in the proof of Theorem~2 of Beyersmann et al. (2013), we have
 $$
  \max_{i\leq 2n} \sqrt{2n} \| (Z_{2n;i}(t_1), \dots, Z_{2n;i}(t_k)) \| = o_P(1),
 $$
 where $\|\cdot\|$ denotes the euclidean distance. This implies condition~$(4.1)$ in Pauly (2011). Now the calculation of $(4.2)$ in Pauly (2011) finishes the proof of the finite dimensional convergence:
 The matrix
$$
  \sum_{i=1}^{2n} 2n \left((Z_{2n;i}(t_j))_j - (\bar Z_{2n}(t_j))_j\right) 
  \left((Z_{2n;i}(t_\ell))_\ell - (\bar Z_{2n}(t_\ell))_\ell\right)^T
 $$
 has the entries
 \begin{eqnarray}
  \label{eq:conv_covariance}
   && 2n \sum_{i=1}^{n} [X_{n;i}(t_j) X_{n;i}(t_\ell) + Y_{n;i}(t_j) Y_{n;i}(t_\ell)]\nonumber\\
   && - \sum_{i=1}^{n} [X_{n;i}(t_j) + Y_{n;i}(t_j)] \sum_{i=1}^{n} [X_{n;i}(t_\ell) +  Y_{n;i}(t_\ell)].
 \end{eqnarray}
 Similarly as in Beyersmann et al. (2013) the first sum converges to $2 \zeta(t_j,t_l)$ in probability. Moreover, each factor of the second sum has the local martingale representation
 \begin{eqnarray}
  && \sum_{i=1}^{n} [X_{n;i}(s) + Y_{n;i}(s)] = \int_0^s \frac{\hat S_2(u-)}{Y(u)} J(u) \, \d M_{1}(u)\nonumber\\
    \label{eq:fidi_conv_1}&& + \int_0^s \frac{\hat F_1(u-)}{Y(u)} J(u) \, \d M_{2}(u) 
    - \hat F_1(s) \int_0^s \frac{J(u)}{Y(u)}  \, \d M_{\cdot}(u) \\
   && + \int_0^s \{ \hat S_2(u-) - \hat F_1(s) \} J(u) \alpha_{1}(u) +
    \{ \hat F_1(u-) - \hat F_1(s) \} J(u) \alpha_{2}(u) \, \d u,\nonumber
  \end{eqnarray}
  where $M_{\cdot} = M_{1} + M_{2} = \sum_{j=1}^2 \sum_{i=1}^n ( N_{j;i} + \int_0^\cdot \alpha_jY_{i} \, \d \leb )$
  is the Doob-Meyer local martingale representation of the counting process $N_1 + N_2$. 
  Note that each of the three first integrals in~\eqref{eq:fidi_conv_1} also is a local square integrable martingale by Theorem~II.3.1 of Andersen et al. (1993). 
  By Rebolledo's martingale limit theorem it is easy to see that each local martingale in~\eqref{eq:fidi_conv_1} converges to zero in probability:
  Consider, for instance,
  \begin{eqnarray*}
   \left< \int_0^\cdot \frac{\hat S_2}{Y} J \, \d M_{1} \right>(s)
    = \int_0^s \frac{ S^{2}_2 }{Y^2} J  \, \d \left< M_{1} \right>
    = \int_0^s \frac{ S_2^{2} }{Y} J \, \d A_{1} 
    \leq \int_0^s \frac{ J }{Y}  \, \d A_{1} 
   \prob 0
    \end{eqnarray*}
  by Condition~\eqref{eq:mainass}, where we have implicitely used the notation of Andersen et al. (1993). 
 A similar result holds for the other local martingales. The remaining integrals, however, converge to 
  $$
   \int_0^s \{ S_2(u) - F_1(s) \} \alpha_{1}(u) \, \d u \;
   \text{ and } \; \int_0^s \{ F_1(u) - F_1(s) \} \alpha_{2}(u) \, \d u
  $$
  in probability by the uniform consistency of the Aalen-Johansen estimator and Condition~\eqref{eq:mainass}, respectively.
  This shows $(4.2)$ in Pauly (2011) and thus the desired finite dimensional convergence. \\
  It remains to prove the conditional tightness of the process. To this end we apply Theorem~13.5 in Billingsley (1999) and rewrite
  $$
    \widehat{W}_n^*(u) = \widehat{W}_n^*((Z_{2n;i})_i)(u) = \sqrt{2n} \sum_{i=1}^{2n} \left( w_{2n;i} - \bar w_{2n} \right)
    Z_{2n;i}(u).
  $$
   Let $0 \leq r \leq s \leq u \leq t$ and $\beta = 1$. Then, by the measurability of $\vZ_{2n}$ and their independence of $\vw_{2n}$, it follows that
  \begin{align}
  \begin{split}
  \label{ref:cond_moment}
   & E \left[ (\widehat{W}_n^*(u) - \widehat{W}_n^*(s))^2 (\widehat{W}_n^*(s) - \widehat{W}_n^*(r))^2 \left| \vZ_{2n} \right. \right]
 \\
   & = 4 n^2 \sum_{i_1,i_2,j_1,j_2=1}^{2n} \Big( \prod_{k=1,2} (Z_{2n;i_k}(u) - Z_{2n;i_k}(s)) (Z_{2n;j_k}(s) - Z_{2n;j_k}(r)) \Big) \\
   & \quad  \times E \Big[ \prod_{\ell=1}^2 (w_{2n;i_\ell} - \bar w_{2n})(w_{2n;j_\ell} - \bar w_{2n}) \Big] \\
   & \leq C_1 D_1 \big| E [ (w_{2n;1} - \bar w_{2n})^4 ] \big|
      + C_2 D_2 \big| E [ (w_{2n;1} - \bar w_{2n})^3 (w_{2n;2} - \bar w_{2n}) ] \big| \\
   & \quad   + C_3 D_3 \big| E [ (w_{2n;1} - \bar w_{2n})^2 (w_{2n;2} - \bar w_{2n})^2 ] \big| \\
   & \quad   + C_4 D_4 \big| E [ (w_{2n;1} - \bar w_{2n})^2 (w_{2n;2} - \bar w_{2n}) (w_{2n;3} - \bar w_{2n})] \big| \\
   & \quad   + C_5 D_5 \Big| E \Big[ \prod_{i=1}^4 (w_{2n;i} - \bar w_{2n}) \Big] \Big| = \sum_{k=1}^5 C_k D_k E_k, \\
  \end{split}
  \end{align}
  where $C_k, k=1,\dots,5$, counts the number of possible index values each leading to the same expected value.
  For example, $C_3 = 3$ due to the index combinations $i_1 = i_2 \neq j_1 = j_2$, $i_1 = j_1 \neq i_2 = j_2$ and $i_1 = j_2 \neq j_1 = i_2$.
  The $D_k$ are defined as
  \begin{align*}
   D_k = \max_{\substack{(x_\ell,y_\ell) \in \{ (r,s), (s,u) \}, \\ \ell= 1, \dots, 4}}
    4n^2 \sum \prod_{\ell = 1}^4 | Z_{2n;i_\ell}(y_\ell) - Z_{2n;i_\ell}(x_\ell) |,
  \end{align*}
  where the sum runs over all indices $i_1,i_2,i_3,i_4$ that yield the expected value $E_k$.
  Each case $k=1,\dots,5$ is treated separately:
  Recall that each $Z_{2n;i}$ is represented by a one-jump process $N_{1;i}$ or $N_{2;i}$ so that
  \begin{align*}
  	D_1 \leq n^2 \sum_{i=1}^{2n} \int_0^u \frac{J \d ( N_{1;i} + N_{2;i} )}{Y^4} \cdot \mathcal{O}(1)
  	= \frac{1}{n} \int_0^u \frac{J}{(Y/n)^3} \d (\hat A_1 + \hat A_2) \cdot \mathcal{O}(1)
  \end{align*}
  which tends to zero in probability by Lemma~\ref{lem:Variances}.
  Condition~\eqref{G5} yields the negligibility of $C_1 D_1 E_1$.\\
  For treating $k=2$ first note that, by the Cauchy-Schwarz inequality,
  $$
    \sum_{i=1}^{2n} |Z_{2n;i}(y) - Z_{2n;i}(x)| \leq \left( 2n \sum_{i=1}^{2n} (Z_{2n;i}(y) - Z_{2n;i}(x))^2  \right)^{1/2}
  $$
  for all $(x,y) \in \{ (r,s), (s,u) \}$.
  It follows that
  \begin{align*}
  	D_2 & \leq \max_{\substack{(x,y) \in \{ (r,s), (s,u) \}}} 4n^2 \sum_{i=1}^{2n} |Z_{2n;i}(y) - Z_{2n;i}(x)|^3 \Big( 2n \sum_{j=1}^{2n} (Z_{2n;j}(y) - Z_{2n;j}(x))^2 \Big)^{1/2} \\
  	& \leq \max_{\substack{(x,y) \in \{ (r,s), (s,u) \}}} \Big( n \sum_{i=1}^{2n} (Z_{2n;i}(y) - Z_{2n;i}(x))^2 \Big)^{3/2} \cdot \mathcal{O}_P(1),
  \end{align*}
  where, by Assumption~\eqref{eq:mainass} and the involved $(Y/n)^{-1}$ in the integrand, the asymptotic boundedness of $\max_i n|Z_{2n;i}(y) - Z_{2n;i}(x)|$ in probability
  yields the last inequality.
  Applying the H\"older$(p,q)$-inequality with $p=3/4, q=1/4$ to the expectation $E_2$, we arrive at an upper bound for $C_2 D_2 E_2$.
  Now Conditions~\eqref{G5} -- \eqref{G7} and straightforward applications of the Cauchy-Schwarz inequality as above imply
  \begin{align*}
    \sum_{k=3}^5 C_k D_k E_k \leq \Big( n \sum_{i=1}^{2n} (Z_{2n;i}(y) - Z_{2n;i}(x))^2 \Big)^2 \cdot \mathcal{O}(1) \\
    \leq \Big( n \sum_{i=1}^{2n} (Z_{2n;i}(y) - Z_{2n;i}(x))^2 \Big)^{3/2} \cdot \mathcal{O}_P(1),
  \end{align*} 
  where $\mathcal{O}_P(1)$ can be chosen independently of $r,s,u$.
  Thus, we have found a common upper bound for $C_k D_k E_k, k=1,\dots,5,$ that equals $\mathcal{O}_P(1)$ times
  \begin{align*}
  	h_n^{3/2}(x,y) := \Big( n \sum_{i=1}^{2n} (Z_{2n;i}(y) - Z_{2n;i}(x))^2 \Big)^{3/2}.
  \end{align*}
  If, for example, $(x,y) = (r,s)$, then $h_n(r,s)$ equals
  $$
    n \sum_{i=1}^{n} ( X_{2n;i}(s) - X_{2n;i}(r) )^2 
      + n \sum_{i=1}^{n} ( Y_{2n;i}(s) - Y_{2n;i}(r) )^2. \\
  $$
Due to similarity, we only consider the first term. Since $N_{1;i}, 1\leq i\leq n,$ are all one-jump processes, this term is equal to
  \begin{eqnarray*}
       && n \sum_{i=1}^{n} \left( \int_r^s \frac{(\hat S_2 - \hat F_1(r)) J \ d N_{1;i}}{Y} - 
      \left( \hat F_1(s) - \hat F_1(r) \right) \int_0^s \frac{J \ d N_{1;i}}{Y} \right)^2 \\
     && \leq 2 n \sum_{i=1}^{n} \left\{ \int_r^s \frac{(\hat S_2 - \hat F_1(r))^2 J \ d N_{1;i}}{Y^2}  
     + \left( \hat F_1(s) - \hat F_1(r) \right)^2 \int_0^s \frac{J \ d N_{1;i}}{Y^2} \right\} \\
     && \leq 2 \left\{ n \left(\hat \sigma_1^2(s) - \hat \sigma_1^2(r) \right) 
      + \left( \hat F_1(s) - \hat F_1(r) \right)^2 n \hat \sigma_1^2(s) \right\},
  \end{eqnarray*}
  where  the left-continuity of all integrands should be kept in mind and $\hat \sigma_1^2(u) = \int_0^u J/Y^2 \d N_{1}$ as in Beyersmann et al. (2013).
  Now $(a - b)^2 \leq a^2 - b^2$ for all $0 \leq b \leq a$ yields the upper bound
  $$
   2 \left\{ n \left(\hat \sigma_1^2(u) - \hat \sigma_1^2(r) \right) 
   + \left( \hat F_1^2(u) - \hat F_1^2(r) \right) n \hat \sigma_1^2(t) \right\}
  $$
  which, by Theorems~IV.1.2 and~IV.4.1 in Andersen et al. (1993), converges uniformly in $r,u \in [0,t]$ to
$$
   2 \left\{ \left(\sigma_1^2(u) - \sigma_1^2(r) \right)  + \left( F_1^2(u) - F_1^2(r) \right) 
   \sigma_1^2(t) \right\},
$$
where $\sigma_j^2(s) = \int_0^s \alpha_j(v)/y(v) dv$ for $j=1,2$, see Equation $(4.1.11)$ in Andersen et al. (1993). 
  Similarly, the convergence of the second sum holds with $\sigma_2^2$ instead of $\sigma_1^2$. We can now finish the proof as in Beyersmann et al. (2013) by the subsequence principle for convergence in probability:
  For each subsequence there exists a further subsequence such that for $P$ a.e. $\omega \in \Omega$ there exists an $\varepsilon > 0$ such that \eqref{ref:cond_moment} 
is less than or equal to $C (H(u) - H(r))^{3/2}$ for large $n \geq n_0$ and a constant $C > 0$.
Note that $\varepsilon, n_0$ and $C$ are independent of $r,s,u \in [0,t]$.
  Here the non-decreasing, continuous function $H$ is given by 
  $$
    H(v) = \left( \sigma_1^2(v) + \sigma_2^2(v) \right)
      + F_1^2(v) \left( \sigma_1^2(t) + \sigma_2^2(t) \right)  + \varepsilon v.
  $$
  Hence the conditional tightness follows from Theorem 13.5 in Billingsley (1999) pointwise along subsequences which in turn implies the assertion of this theorem. \hfill$\Box$

%
{\it Proof of Theorem~\ref{theo: tests}}
As already outlined above the convergences $T_n \oDo T\sim N(0,\sigma_{\zeta}^2)$ and $V_n^2 \oPo \sigma_{\zeta}^2$ (see Lemma~\ref{lem:Variances} below) hold under $H_=$.
Moreover, $\sigma_{\zeta}^2>0$ holds since it is assumed that $\alpha_1^{(k)}>0$ on a set with positive Lebesgue-$\leb_{|[t_1,t_2]}$ measure for at least one choice of $k=1,2$. 
Hence $T_{n,stud}$ is asymptotically standard normal by Slutzky's Lemma. 
In addition, since $\sigma_{\zeta}^2>0$ even holds  for $F_1^{(1)}\neq F_1^{(2)}$, we have that $T_{n,stud}\oPo  \infty \mathbf{1}_{K_{\gneqq}} - \infty \mathbf{1}_{H_{\lneqq}},$   where $H_\lneqq : \{ F_1^{(1)}\lneqq F_1^{(2)} \text {on } [t_1,t_2]\}$. 
Altogether this proves the consistency and asymptotic exactness of $\varphi_n$ under $H_=$. It remains to investigate the conditional asymptotic behaviour of $T_{n,stud}^*$. 
To this end, Theorem~\ref{theo:exchCLT} together with Example~\ref{ex.weights} and the continuous mapping theorem show that 
the conditional distribution of $T_n^*$ given the data is asymptotically $N(0,\sigma_{\tilde{\zeta}}^2)$-distributed with 
$$
  \sigma_{\tilde{\zeta}}^2 = \int_{t_1}^{t_2} \int_{t_1}^{t_2} \rho(r) \Big(\big(p_2\zeta^{(1)} + p_1\zeta^{(2)}\big) (r,s) - \frac{p_1p_2}{2} \big(\xi^{(1)} - \xi^{(2)}\big)(r) \big(\xi^{(1)}-\xi^{(2)}\big)(s)   \Big)\rho(s) \, \d r \d s.
$$
Note, that the continuous mapping theorem is indeed applicable since there exist versions of $U^{(k)}$ and $U^{*\, (k)}, k=1,2,$ with a.s. continuous sample paths. 
Moreover, it is proven in Lemma~\ref{lem:Variances} that $V_n^{*\,2}$ converges in probability to $\sigma_{\tilde{\zeta}}^2$ which is positive by assumption. 
Hence it follows that the conditional distribution of $T_{n,stud}^*$ given the data is asymptotically standard normal. As above this proves consistency and asymptotic exactness under $H_=$ of $\varphi_n^*$. 
Finally, the asymptotic equivalence of both tests follows from Lemma~1 in Janssen and Pauls (2003).
\hfill$\Box$
%

%
\begin{lemma}\label{lem:Variances}
 (a) With the notation of Section~\ref{sec: model} suppose that Condition~\eqref{eq:mainass} holds. Then for $t<\tau$ and for all $r < \ell - 1$ and $j=1,2,$ the stochastic process 
 \begin{align*}
   \Big( \hat{\sigma}(s) := n^r \sum_{i=1}^n \int_0^s h(u) \frac{J(u)}{Y^\ell(u)} \d N_{j;i}(u) \Big)_{s \in [0,t]} 
 \end{align*}
converges to zero on $D[0,t]$ in probability if the left-continuous function $h$ is bounded by a constant $C > 0$.\\[0.25ex]
 (b) Under the assumptions of Theorem~\ref{theo: tests} the variance estimators $V_n^{2}$ and $V_n^{*\,2}$ are consistent estimates for $\sigma_{\zeta}^2$ and $\sigma_{\tilde{\zeta}}^2$, respectively.
\end{lemma}

{\it Proof of Lemma~\ref{lem:Variances}} 
(a) Clearly, $\hat{\sigma}$ is bounded by a process with Doob-Meyer decomposition
 \begin{align*}
  | \hat{\sigma}(s) | \leq C n^r \sum_{i=1}^n \int_0^s \frac{J}{Y^\ell} \d N_{j;i}
    = C n^r \int_0^s \frac{J}{Y^\ell} \d M_j + C n^r \int_0^s \frac{\alpha J }{Y^{\ell-1}} \d\leb,
 \end{align*}
 where $M_j = \sum_{i=1}^n M_{j;i}$ are locally square integrable martingales. 
 The local martingale in the above decomposition has the predictable covariation process
 \begin{align*}
  \left< C n^r \int_0^\cdot \frac{J}{Y^\ell} \d M_j \right>(s)
  = C^2 n^{2r} \int_0^s \frac{\alpha J }{Y^{2\ell-1}} \d\leb .
 \end{align*}
 Both this expression and $n^r \int_0^s \alpha J /Y^{\ell-1} \d \leb$ converge to zero in probability as $n \rightarrow \infty$ if $r < \ell-1$.
 Eventually, Rebolledo's Theorem yields the desired convergence on $D[0,t]$.\\
(b) Note first that the processes $\hat \zeta_n$ and $\hat \xi_n :=\sqrt{\frac{n_1n_2}{2n^2}}\sum_{i=1}^{2n} Z_{2n;i}$ can be decomposed into several additive, monotonic functions on $[t_1,t_2]^2$
each of which converges (pointwise on $[t_1,t_2]^2$) towards its real, unknown, monotonic and continuous counterpart in probability as $n \rightarrow \infty$.
This is due to the consistency of the Aalen-Johansen estimator for CIFs as well as a similar argument as in Beyersmann et al. (2013).
A simple Polya-type argument now shows that such monotonic process estimators even converge uniformly on $[t_1,t_2]^2$ in probability
which implies the convergence of the weighted integrals over $\hat \zeta_n$ and $\hat \xi_n(r) \hat \xi_n(s)$, in particular the convergence of $V_n^2$ in probability.\\
We now continue to show the consistency of $V_n^{*\, 2}$ and start by proving that 
 \begin{align}
 \label{eq:int_cov_est_second_moment_0}
   \E\left[ \left. \left( \int_{t_1}^{t_2} \int_{t_1}^{t_2} \rho(r) \left( \zeta^*(r,s) 
   - \hat \zeta(r,s) \right) \rho(s) \d r \d s \right)^2 \right| \vZ_{2n} \right]
  \end{align}
is negligible. Recall, that  $Z_{2n;i}$ are defined as integrals with respect to counting processes. 
We now pool each quantity in a canonical way by merging the indices $k$ and $i$, i.e. $(v_{2n;\ell})_\ell = (v_{2n;i}^{(k)})_{i,k},$
$(N_{\ell})_\ell = (N_{1;i}^{(k)} + N_{2;i}^{(k)})_{i,k}$ and similarly for $J$ and $Y$. 
Then, after changing the order of integration to $\d r \d s \d N_{j;i}^{(k)}$, we see that~\eqref{eq:int_cov_est_second_moment_0} is bounded from above by
  \begin{align}
  \label{eq:int_cov_est_second_moment_1}
   \left( \frac{n_1 n_2}{n} \right)^2 \sum_{\ell_1,\ell_2}^{2n}
   \int_0^{t_2} \frac{h_{\ell_1} J_{\ell_1}}{Y_{\ell_1}^2} \d N_{\ell_1}  \int_0^{t_2} \frac{h_{\ell_2} J_{\ell_2}}{Y_{\ell_2}^2} \d N_{\ell_2}
   \left|\E[(v_{2n;\ell_1} - 1) (v_{2n;\ell_2} - 1)]\right|,
  \end{align}
where 
  \begin{align*}
   h_{l_k}(u) := & \iint\limits_{[u \wedge t_1, t_2]^2}
     \textbf{1}(l_k \leq n_1)
    (\hat S_2^{(1)}(u) - \hat F_1^{(1)}(r))(\hat S_2^{(1)}(u) - \hat F_1^{(1)}(s)) \\
    & + \textbf{1}(n_1 < l_k \leq n)
    (\hat F_1^{(1)}(u) - \hat F_1^{(1)}(r))(\hat F_1^{(1)}(u) - \hat F_1^{(1)}(s)) \\
    & + \textbf{1}(n < l_k \leq n + n_1)
    (\hat S_2^{(2)}(u) - \hat F_1^{(2)}(r))(\hat S_2^{(2)}(u) - \hat F_1^{(2)}(s)) \\
    & + \textbf{1}(n + n_1 < l_k)
    (\hat F_1^{(2)}(u) - \hat F_1^{(2)}(r))(\hat F_1^{(2)}(u) - \hat F_1^{(2)}(s)) \d r \d s
\end{align*}
are bounded functions. Straightforward calculations show that 
  \begin{align*}
    C:= \limsup_{n \rightarrow \infty} \left| \E[(v_{2n;\ell_1}-1)(v_{2n;\ell_2}-1)] \right| (n \b 1(\ell_1\neq \ell_2) + \b 1(\ell_1=\ell_2)) <\infty
  \end{align*}
holds for both choices of $v_{2n;\ell}$ (i.e. in Efron's or the wild bootstrap case).
Hence, for large $n$, the absolute value of~\eqref{eq:int_cov_est_second_moment_1} has the upper bound
  \begin{align*}
   (C+1) p_1^2 p_2^2 \left\{ n^2 \sum_{\ell=1}^{2n} \int_0^{t_2} \frac{h_\ell^2 J_\ell}{Y_\ell^4} \d N_\ell
    + \left( n^{1/2} \sum_{\ell=1}^{2n} \int_0^{t_2} \frac{|h_\ell| J_\ell}{Y_\ell^2} \d N_\ell \right)^2 \right\}.
  \end{align*}
  Part (a) now yields the convergence of 
\begin{equation}\label{proof: int}
  \int_{t_1}^{t_2} \int_{t_1}^{t_2} \rho(r) (\hat \zeta_{n} - \zeta_{n}^*) (r,s) \rho (s) \d r \d s
\end{equation}
to zero in probability given the data. In the same way it can be shown that the remaining integral with $(\hat \zeta_{n} - \zeta_{n}^*) (r,s)$  replaced by  $\xi_{n}^*(r,s)-\hat\xi_n(r)\hat\xi_n(s)$ in \eqref{proof: int} also converges 
to zero in probability given the data which completes the proof.\hfill$\Box$\\

Finally, we give the examples mentioned in Remark~\ref{remark:theo}(a) and prove that they fulfill the assumptions of Theorem~\ref{theo:exchCLT}. 
The extensions to the two-sample case as mentioned in Section~\ref{sec: Test} are straightforward.

\begin{example}
 \label{ex.weights}
For the following resampling weights the convergence \eqref{eq:weakexchnBSconv} from Theorem~\ref{theo:exchCLT} is fulfilled.\\[1ex]
(a) Let $(m_{2n;1},\dots,m_{2n;2n})$ be a multinomially distributed random vector with sample size $2n=\sum_{i=1}^{2n} m_{2n;i}$ and equal selection probability $1/2n$. 
Then {\it  Efron's classical bootstrap} weights 
\begin{equation}\label{eq:mnBS}
  w_{2n;i} = m_{2n;i} - 1,\quad 1\leq i\leq 2n,
\end{equation}
are covered by our approach.\\[0.5ex]
(b) Let $G_{2n;i}$ be row-wise i.i.d. weights with $\limsup_{n} \E(G_{2n;1}^4) <\infty$ as well as $\E(G_{2n;1})=0, Var(G_{2n;1})=1$. Then the {\it wild bootstrap} weights 
\begin{equation}\label{eq:wildBS}
 w_{2n;i} = G_{2n;i},\quad 1\leq i\leq 2n,
\end{equation}
fulfill the Conditions \eqref{G2} -- \eqref{G7} provided that $G_{2n;1}\dist Z$.\\[0.5ex]
(c) As special example the choice $G_{2n;i}=G_i-1$ for i.i.d. $Poi(1)-$distributed random variables $G_1,\dots,G_{2n}$ yields the so called {\it Poisson-bootstrap} which may be interpreted as drawing $N=\sum_{i=1}^{2n} G_i$ times with replacement from $\vZ_{2n}(\cdot)$. 
Moreover, the choice $G_{2n;i}=G_i'$ for $G_i'\stackrel{\text{i.i.d.}}{\sim}N(0,1)$ corresponds to {\it Lin's resampling technique}. \\[0.75ex]
(d) Let $\eta_i>0, 1 \leq i \leq 2n,$ be positive i.i.d. random variables with $E(\eta_1)=\mu_\eta, \sigma_\eta^2=Var(\eta_{1})$ and finite fourth moment. 
Then the limit Theorem~\eqref{eq:weakexchnBSconv} holds for the {\it i.i.d. weighted bootstrap} weights 
$w_{2n;i}= C_\eta^{-1} ({\eta_{i}}/{\overline{\eta}_{2n}} - 1),$
where $C_\eta^2 = \sigma_\eta^2/\mu_\eta^2,$ and $\overline{\eta}_{2n} = \sum_{i=1}^{2n} \eta_i/2n$.\\[0.75ex]
(e) Rubin's {\it Bayesian bootstrap} is achieved by letting $\eta_i\stackrel{\text{i.i.d.}}{\sim}Exp(1)$ in (d) with $C_\eta=1$.

\end{example}

{\it Proof of Example~\ref{ex.weights}.} 
We first show that the weights given in (a)--(c) fulfill the Conditions \eqref{G2} -- \eqref{G7}. 
Since thereof part (a) is the most difficult to prove, we only consider this part and leave the others as an exercise. 
Moreover, we only show that Condition~\eqref{G7} holds, since~\eqref{G5} and~\eqref{G6} can be shown similarly and the prove for \eqref{G2} -- \eqref{G4} can be found in Janssen (2005) and Pauly (2009). 
Let $n \geq 2$, then we start with
 \begin{align*}
  & \E \Big( \prod_{i=1}^4 (m_{2n;i} - \bar m_{2n}) \Big)
  = \E \Big( \prod_{i=1}^4 (m_{2n;i} - 1) \Big) \\
  & = \E \Big( \prod_{i=1}^4 m_{2n;i} \Big)
    - 4 \E \Big( \prod_{i=1}^3 m_{2n;i} \Big)
    + 6 \E \Big( \prod_{i=1}^2 m_{2n;i} \Big) 
- 4 \E \Big( m_{2n;1} \Big) + 1
  \end{align*}
 where each single expectation is further calculated with the help of the moment generating function of $(m_{2n;i})_i$ or by consulting the monograph of Johnson et al. (1997) \nocite{johnsonkotz_1997}. 
 Thus, we have
 \begin{align*}
    & \E \Big( \prod_{i=1}^4 m_{2n;i} \Big) = \frac{2n(2n-1)(2n-2)(2n-3)}{16n^4}, \\
    & \E \Big( \prod_{i=1}^3 m_{2n;i} \Big) = \frac{2n(2n-1)(2n-2)}{8n^3}
 \end{align*}
 and $\E[m_{2n;1} m_{2n;2}] = cov(m_{2n;1}, m_{2n;2}) + \E[m_{2n;1}]^2 = - 2n \frac{1}{4n^2} + 1 = 1 - \frac{1}{2n}$
 so that the initial expectation finally equals
 \begin{align*}
  & \frac{2n(2n-1)(2n-2)(2n-3)}{16n^4} - 4 \frac{2n(2n-1)(2n-2)}{8n^3} + 6 \Big(1 - \frac{1}{2n}\Big) - 3 \\
  & = \frac{3}{4n^2} - \frac{3}{4n^3} \in \mathcal{O}(n^{-2}).
 \end{align*} 
Hence (a) follows. Part (b) can be shown in the same way and (c) is only a special example of (b). 
We will now prove (d) with the help of (b). To this end we rewrite $\widehat{W}_n^*$ as
\begin{eqnarray*}
 \widehat{W}_n^* &=& C_{\eta} \sqrt{2n} \sum_{i=1}^{2n} \frac{\eta_i}{\overline{\eta}_{2n}}(Z_{2n;i} - \overline{Z}_{2n}) 
\ = \ \frac{C_{\eta}\sigma_\eta}{\overline{\eta}_{2n}} \sqrt{2n} \sum_{i=1}^{2n} \frac{(\eta_i - \mu_\eta)}{\sigma_\eta} (Z_{2n;i} - \overline{Z}_{2n}),
\end{eqnarray*}
where we have utilized in the first and last equality the identity $\sum_i(Z_{2n;i} - \overline{Z}_{2n}) =0$. Here the first factor ${C_{\eta}\sigma_\eta}/{\overline{\eta}_{2n}}$ on the right hand side converges to $1$ 
almost surely by the law of large numbers and the second factor is a wild bootstrap version \eqref{eq: exchBSAJ} of the Aalen-Johansen estimator in the weights $G_i = {(\eta_i - \mu_\eta)}/{\sigma_\eta}$. 
Hence the assertion is a consequence of Slutzky's Lemma and part (b). Part (e) is only a special example of (d).
\hfill$\Box$
%
%

\bibliography{literatur}
\bibliographystyle{plain}

\end{document}